\documentstyle[12pt]{article}

\newcommand{\be}{\begin{equation}}
      \newcommand{\ee}{\end{equation}}
      \newcommand{\ba}{\begin{eqnarray}}
       \newcommand{\ea}{\end{eqnarray}}
\newcommand{\ban}{\begin{eqnarray*}}
       \newcommand{\ean}{\end{eqnarray*}}

\newcommand{\lp}{\langle}
\newcommand{\rp}{\rangle}
\newcommand{\ra}{\rightarrow}

 \renewcommand{\o}[2]{\frac{#1}{#2}}
\newcommand{\hf}{\o{1}{2}}

 \newcommand{\qed}{\hspace*{\fill}\rule{3mm}{3mm}\quad}
 \newcommand{\Pf}{\noindent {\em Proof.} }

\newcommand{\Rk}{\noindent {\em Remark} }

\newcommand{\sect}[1]{\section{#1} \setcounter{equation}{0}}

\newtheorem{theo}{Theorem}[section]

\begin{document}
\newtheorem{lem}[theo]{Lemma}
\newtheorem{prop}[theo]{Proposition}  
\newtheorem{coro}[theo]{Corollary}

\title{Ricci Curvature and Betti Numbers}
\author{Guofang Wei\thanks {Supported in part by NSF DMS-9409166}\\
Department of Mathematics \\
 University of California \\
Santa Barbara, CA 93106 \\
 wei@math.ucsb.edu}
\date{}
\maketitle

\begin{abstract}
We derive a uniform bound for the total betti number of a closed manifold in terms of a Ricci curvature lower bound, a conjugate radius lower bound and a diameter upper bound. The result is based on an angle version of Toponogov comparison estimate for small triangles in a complete manifold with a Ricci curvature lower bound.  We also give a uniform estimate on the generators of the fundamental group and prove a fibration theorem in this setting.
\end{abstract}

\newcommand{\inj}{\mbox{inj}}
\newcommand{\vol}{\mbox{vol}}
\newcommand{\diam}{\mbox{diam}}
\newcommand{\Ric}{\mbox{Ric}}
\newcommand{\Iso}{\mbox{Iso}}
\newcommand{\Hess}{\mbox{Hess}} 
\newcommand{\divg}{\mbox{div}}
\sect{Introduction}
%\pagestyle{myheadings}
%\markright{A Comparison-Estimate of Toponogov Type for Ricci Curvature}

Let $M$ be a compact Riemannian manifold of dimension $n$. In \cite{g} Gromov proved the following celebrated result for manifolds with lower sectional curvature bound.
\begin{theo}[Gromov] \label{gtheo}
Given $n, D > 0,\ H$, and a field $F$, if
\[
 K_M \geq H, \ \ \diam (M^n) \leq D, 
\]
then 
\be   \label{bg}
\sum_i b^i(M^n; F) \leq C(n, HD^2).
\ee
\end{theo}

For manifolds with just Ricci curvature bounded below Gromov showed that (see \cite{glp}) if $\Ric_M \geq H, \ \diam (M^n) \leq D$, then $b_1(M; R) \leq C(n, HD^2)$. For higher betti numbers this estimate does not hold anymore (see e.g. \cite{sy}). Recently Perelman \cite{p1} showed that the estimate (\ref{bg}) is not valid even with an additional lower bound on volume. (In this case the  manifolds can not collapse.)

For the class of manifolds $M^n$ satisfying 
\be \label{class}
\Ric_M \geq (n-1)H, \ \mbox{conjugate radius} \ \geq r_0,\ \diam (M) \leq D.
\ee
we have the following estimate for betti numbers.
\begin{theo}[Uniform betti number estimate] \label{bn}
For the class of  manifolds $M^n$ satisfying (\ref{class}), we have
\[
\sum_i b^i (M^n) \leq C(n,H,r_0,D).
\]
\end{theo}

\Rk The class of manifolds satisfying (\ref{class}) could have infinite many different homotopy types, e.g. the space forms.

We also have
\begin{theo} \label{fg}
For the class of  manifold $M^n$ satisfying (\ref{class}), the fundamental group can be generated by $s \leq s(n,H,r_0,D)$ elements, where $s(n,H,r_0,D)$ is a constant depends only on $n,H,r_0,D$.
\end{theo}

\begin{theo} \label{fib}
Given $n,\ i_0 >0,\ r_0 >0$, there exists a positive number $\epsilon (n,i_0,r_0)$ such that if complete manifolds $M^n,\ N^m$ satisfying Ric$_M \geq -1$, conj$_M \geq r_0$, $|K_N| \leq 1$, inj$_N \geq i_0$ and $d_H (M,N) \leq \epsilon (n,i_0,r_0)$, then there exsits a fibration $f: M \ra N$.
\end{theo}

\begin{prop} \label{sph}
Let $M^n$ be a closed  manifold with $\Ric \geq (n-1)$, conjugate radius $\geq r_0$. Let $p,q \in M$ be such that $|pq| =$ diam$_M$. Then for any $\delta > 0$ there is a constant $\epsilon = \epsilon (n,\delta,r_0) > 0$ such that if diam $\geq \pi - \epsilon$, then any $x \in M -\{ B_\delta (p) \cup B_\delta (q) \}$ is a regular point for $p$ as well as for $q$.
\end{prop}

 Theorem~\ref{bn}, \ref{fib} is not true without the conjugate radius lower bound (see \cite{sy,a2}). But the results all hold with sectional curvature lower bound instead of conjugate radius lower bound (see \cite{g}, \cite{y}, \cite{gp}). In this case Theorem~\ref{fib} can be strengthened to that $f$ is almost a Riemannian submersion and the fibre is a manifold of almost nonnegative  sectional curvature (see \cite{y} for detail, see also \cite{f, cfg}), Proposition~\ref{sph} can be strengthened to that it is true for $\delta = 0$,  proving that $M^n$ is a twisted sphere (see \cite{p2}).
So far the only natural way to obtain conjugate radius lower bound is via sectional curvature upper bound. But one can easily construct examples of manifolds satisfying (\ref{class}) where the sectional curvatures are not uniformly bounded.

The class of manifolds satisfying (\ref{class}) could collapse, e.g. the space forms.
In the noncollapsing case (i.e. with an additional lower bound on the volume) this class is well understood. In fact by \cite{cgt} the injectivity radius can be bounded below by Ricci curvature, conjugate radius and volume lower bounds. (This is pointed out to me by Peter Petersen.) Combining this with \cite{ac} one knows that the class of manifolds satisfying (\ref{class}) and with volume bounded below are $C^\alpha$ compact.  See \cite{a,w} for other results in the non-collapsing case. Thus our results here can be thought of as a first step in understanding the manifolds with Ricci curvature lower bound which could collapse.

We refer to \cite{co, cc} for some recent significant developments for Ricci curvature.

The essential tool in proving Theorem~\ref{bn}, \ref{fg} and Proposition~\ref{sph} is an angle version Toponogov comparison estimate  for Ricci curvature which we state below. First we introduce some notations.

In the paper, a {\em geodesic triangle} $\{\gamma_0, \gamma_1, \gamma_2\}$ consists of three {\em minimal} geodesics, $\gamma_i$, of length $L[\gamma_i] = l_i$, which satisfy 
\[
\gamma_i (l_i) = \gamma_{i+1} (0) \ \mbox{mod} \ 3 \ (i = 0,1,2).
\]

The {\em angle} at a corner, say $\gamma_0 (0)$, is by definition, $\angle (-\gamma_2'(l_2), \gamma_0'(0))$. The angle opposite $\gamma_i$ will be denoted by $\alpha_i$.

\begin{theo}[Toponogov type comparison-estimate] \label{tce}
For the class of  manifolds $M^n$ satisfying $\Ric_M \geq H$, conj $\geq r_0$ there is a constant $r = r(n,H,r_0) > 0$ such that if $\{\gamma_0, \gamma_1, \gamma_2\}$ is a geodesic triangle contained in $B_p(r)$ for some $p \in M$, then there is a geodesic triangle $\{\bar{\gamma}_0, \bar{\gamma}_1, \bar{\gamma}_2\}$ in the Euclidean plane with $L[\bar{\gamma}_i] = L[\gamma_i] = l_i$ and 
\be  \label{ae}
 \alpha_i > \frac{18}{19} \bar{\alpha}_i = \frac{18}{19}\cos^{-1} (\frac{l_{i-1}^2 + l_{i+1}^2 - l_i^2}{2l_{i-1}l_{i+1}}).
\ee 
\end{theo}

\Rk In fact, we can show that for any $0 < \mu < 1$, there is a constant $r = r(n,H,r_0,\mu) > 0$ such that $\alpha_i > \mu \bar{\alpha}_i$ for any geodesic triangles inside $B_p(r)$.

We also refer to \cite{dw} for a hinged version of Toponogov type comparison estimate for Ricci curvature.

The basic structure of the proof of Theorem~\ref{bn},  \ref{fg}, \ref{fib} and Proposition~\ref{sph} is the same as in the case with sectional curvature lower bound. In the presence of sectional curvature lower bound we have the powerful  Toponogov comparison theorem. Here we try to replace that by Theorem~\ref{tce}, the Toponogov comparison estimate. There are situations where the rigid structure of the Toponogov comparison theorem is essential. But, as we will see,  the Toponogov comparison estimate suffices for the above results. There is also the problem of going from local to global in applying Theorem~\ref{tce}. That is dealt with differently in each case.

Theorem~\ref{tce} depends essentially on the result of \cite{ac}. To apply the $C^\alpha$-compactness of \cite{ac} we pullback the metric to the tangent space by the exponential map and this is where the conjugate radius comes in. The difficulty with lifting triangles is overcomed by relating the angle comparison radius (see \S 2 below for the definition) on the manifold with the angle comparison radius on the tangent space (with a different base point). We aslo used the regularity result of Calabi-Hartman for geodesics \cite{ch}.

{\bf Acknowledgement}. After we proved the hinged version of the Toponogov comparison estimate [DW]
Shunhui Zhu suggested a rescaling argument (that proves a weaker version).
It turns out that the rescaling argument is the right approach to the 
angle version.  I wish to express my thanks to S. Zhu for that, and other valuable discussions, and  to Xianzhe Dai for reading an earlier version of this paper and constructive suggestions. I would like to thank Rugang Ye for pointing out a gap in the first version of this paper and suggesting the argument to get around it. I am also grateful to M. Anderson, J. Cheeger, D. Gromoll, D. Moore, K. Grove and P. Petersen for helpful discussions and comments.

\sect{Angle comparison for Ricci curvature}
In this section we prove the angle comparison estimate (Theorem~\ref{tce}) following the ideas of \cite{ac}. In order to get the desired $C^\alpha$-convergence we lift everything to the tangent space and this is exactly where the hypothsis on the conjugate radius comes in. We consider balls for which (\ref{ae}) holds for any geodesic triangles inside it and show that we have a uniform positive lower bound for the radius.

\noindent {\em Definition}. Let $M^n$ be a Riemannian manifold and $p \in M$ a fixed point . We define the {\em angle comparison radius} at $p$
 $$r_{ac} (p) = \max \{ r\ | \forall \{\gamma_0, \gamma_1, \gamma_2\} \subset B_p(r),  \alpha_i > \frac{18}{19}\cos^{-1} (\frac{l_{i-1}^2 + l_{i+1}^2 - l_i^2}{2l_{i-1}l_{i+1}}) \}.$$

Note that $r_{ac}(p)$ is strictly positive for any fixed compact smooth Riemannian manifold. As $p$ varies, $r_{ac}(p)$ defines a function of $p$ on $M$. At each point $p \in M$ $r_{ac}(p)$ is the radius of the largest geodesic ball about $p$ on which the angle version of Toponogov estimate (\ref{ae}) holds.

We first prove the following basic property about the angle comparison radius.
\begin{prop}  \label{rc}
If $(M_i, g_i, p_i)$ converges to $(M,g,p)$ in the $C^\alpha$ topology, then
\be \label{rl}
r_{ac}(p) \leq \underline{\lim}_{i \ra \infty} r_{ac} (p_i).
\ee
In particular, if the limit $(M,g,p)$ is the Euclidean space then 
\[
\lim_{i \ra \infty} r_{ac} (p_i) = +\infty.
\]
\end{prop}

To prove Proposition~\ref{rc} one needs the following result of Calabi and Hartman \cite{ch} on the smoothness of isometries.
\begin{theo} [\cite{ch}]
Let $(M,g)$ be a $C^\alpha$ Riemannain manifold with respect to some coordinate. Then its geodesics are $C^{1,\alpha}$ with respect to the same coordinate. Moreover the $C^{1,\alpha}$ norm of the geodesics can be bounded by the $C^\alpha$ norm of the metric.
\end{theo}

An immediate corollary of this result is the following lemma.
\begin{lem}  \label{tc}
If $(M^n_i, g_i, p_i)$ converges to $(M^n,g,p)$ in $C^\alpha$ topology, then for any geodesic triangles $\{\gamma_0^i, \gamma_1^i, \gamma_2^i \}$ in $B_{p_i} (r) \subset M_i$ for some $r > 0$, a subsequence of $\{\gamma_0^i, \gamma_1^i, \gamma_2^i \}$ converges to a geodesic triangle $\{\gamma_0, \gamma_1, \gamma_2 \}$ in $M$ in $C^{1,\alpha'}$ topology, $\alpha' < \alpha$.
\end{lem}

\noindent {\em Proof of Proposition~\ref{rc}} 
If (\ref{rl}) does not hold, then there exsits an $\epsilon > 0$ such that for all $i$ large
\[
r_{ac}(p_i) < r_{ac}(p) -\epsilon.
\]
So there are geodesic triangles $\{\gamma_0^i, \gamma_1^i, \gamma_2^i \}$ in $B_{p_i}(r_{ac}(p)-\frac{\epsilon}{2})$, for which the angle comparison estimate (\ref{ae}) does not hold. By Lemma~\ref{tc} a subsequence of $\{\gamma_0^i, \gamma_1^i, \gamma_2^i \}$ converges to a geodesic triangle $\{\gamma_0, \gamma_1, \gamma_2 \}$ in $B_p(r_{ac}-\frac{\epsilon}{4}) \subset M$. If $\{\gamma_0, \gamma_1, \gamma_2 \}$ is a nontrivial triangle (i.e. at least one of the sides has positive length), then we have a contradiction (since angle converges to angle). In the case $\{\gamma_0, \gamma_1, \gamma_2 \}$ is a point, we rescale the metrics before passing to the limit. That is, we let $q_i$ be one of the edge points of the geodesic triangle $\{\gamma_0^i, \gamma_1^i, \gamma_2^i \}$, $\underline{l}^i = \max \{ l_0^i, l_1^i, l_2^i \}$ and $\bar{g}_i = (\underline{l}^i)^{-2} g_i$. Since $\lim_{i \ra \infty} \underline{l}^i  = 0$ as $i \ra \infty$ $(M_i, \bar{g}_i, q_i)$ converges in $C^\alpha$-topology to the Euclidean space.  Now the geodesic triangles $\{\gamma_0^i, \gamma_1^i, \gamma_2^i \}$ are contained in $B_{q_i}(2)$ in the rescaled metric $\bar{g}_i$, so a subsequence converges to some nontrivial geodesic triangle $\{\bar{\gamma}_0, \bar{\gamma}_1, \bar{\gamma}_2 \} \subset R^n$. But $\{\bar{\gamma}_0, \bar{\gamma}_1, \bar{\gamma}_2 \}$ 
does not satisfy the angle comparison estimate (\ref{ae}), which is also a contradiction.
\qed

We now turn to the proof of the angle comparison estimate.

\noindent {\em Proof of Theorem~\ref{tce}} 
We show that $r_{ac}(p)$ has a uniform lower bound depending only on the bounds
\be  \label{class1}
 \Ric_M \geq H, \ \ \ \mbox{conj} \geq r_0.
\ee
To establish this lower bound on $r_{ac}(p)$, we argue by contradiction. Thus if Theorem~\ref{tce} were false, then there must exist a sequence of Riemannian manifods
$(M_i^n, g_i)$ satisfying (\ref{class1}) but with 
\[
r_i \stackrel{\triangle}{=} r_{ac} (p_i) \ra 0
\]
for some $p_i \in M_i$. Rescale the metric $g_i$ by $r_i^{-2}$, i.e. define metrics $h_i = r_i^{-2} g_i$. Thus $r_{ac} (p_i) = 1$ for $(M_i,h_i)$, while 
\[
\Ric \geq H r_i^2 \ra 0, \ \mbox{conj}_{(M_i,h_i)} \geq r_0/r_i \ra +\infty.
\]
Since $r_{ac} (p_i) = 1$, we have a geodesic triangle $\{\gamma_0^i, \gamma_1^i, \gamma_2^i \} \subset B_{p_i} (2)$ such that one of its angles doesn't satisfy the angle comparison (\ref{ae}). Without loss of generality we assume the angle is at the edge $q_0^i = \gamma_2^i(0)$ and denote it by $\alpha_0^i$.
Let $\tilde{B}_{\tilde{q}_0^i}(r_0/r_i)$ be the ball of radius $r_0/r_i$ in the tangent space $T_{\tilde{q}_0^i}M_i$ with pull back metric $h_i$. We show that 
\be  \label{acru}
\tilde{r}_{ac} (\tilde{q}_0^i) \leq 6.\ee
 To prove this, note that for $i$ large we can lift the geodesics $\gamma_1^i, \ \gamma_2^i$ to minimal geodesics $\tilde{\gamma}_1^i, \ \tilde{\gamma}_2^i$ in $\tilde{B}_{\tilde{q}_0^i}(r_0/r_i)$ with $L(\tilde{\gamma}_1^i) = L(\gamma_1^i), \ L(\tilde{\gamma}_2^i) = L(\gamma_2^i),\ \tilde{\gamma}_1^i (l_1^i)= \tilde{\gamma}_2^i(0) = \tilde{q}_0^i$ and $\tilde{\alpha}_0^i = \angle (-\tilde{\gamma}_1^i (l_1^i), \tilde{\gamma}_2^i(0)) = \alpha_0^i$. Connect $\tilde{\gamma}_1^i (0)$ and  $\tilde{\gamma}_2^i(l_2^i)$ by a minimal geodesic $\tilde{\gamma}_0^i$, then 
\[
L(\tilde{\gamma}_0^i) \geq L(\gamma_0^i).
\] 
(This is suggested to me by Rugang Ye.) This implies the comparison angle for $\tilde{\alpha}_0^i, \ \bar{\tilde{\alpha}}_0^i$, is great or equal to the comparison angle for $\alpha_0^i, \ \bar{\alpha}_0^i$. Therefore
\[
\tilde{\alpha}_0^i = \alpha_0^i \leq \frac{18}{19} \bar{\alpha}_0^i \leq  \frac{18}{19} \bar{\tilde{\alpha}}_0^i,
\]
i.e. the geodesic triangle $\{\tilde{\gamma}_0^i, \tilde{\gamma}_1^i, \tilde{\gamma}_2^i \}$  doesn't satisfy the angle comparison estimate (\ref{ae}), which proves (\ref{acru}).

On the other hand the injectivity radius of $\tilde{B}_{\tilde{q}_0^i}(r_0/r_i)$  is equal to $r_0/r_i$. By \cite{ac} (the version for manifolds with boundary) a subsequence of the manifolds 
$(\tilde{B}_{\tilde{q}_0^i}(r_0/r_i), h_i,p_i)$ converges in the $C^{\alpha'}$ topology, $\alpha' < \alpha$, to a complete (non-compact) $C^\alpha$ Riemannain manifold $(N,h,\tilde{q})$, with $\tilde{q} = \lim \tilde{q}_0^i$. Moreover by \cite[Propositions 1.2, 1.3]{ac} $N$ is isometric to $R^n$, with the canonical flat metric. But $r_{ac} = +\infty$ for $R^n$. By Proposition~\ref{rc} this is contradicting to (\ref{acru}).
\qed

\sect{Uniform betti number estimate}

The proof will use the original ideas of Gromov \cite{g}, but we will follow more closely the beautiful exposition \cite{c}.

To localize the problem, Gromov introduced the {\em content} of a ball,
\[
\mbox{cont} (p,r) = \sum_i \dim  \mbox{Im} (H_i (B_p (r)) \ra H_i (B_p (5r))).\]
Note that if $r > \diam (M)$, then cont $(p,r) = \sum_i b^i(M)$. Thus the basic idea is to estimate cont $(p,r)$ in terms of that of smaller balls (using the Mayer-Vietoris principle).

Now there are two ways of reducing the size of the ball. First by the Mayer-Vietories principle, one has \cite[Corollary 5.7]{c}
\begin{lem} \label{cs}
Let $N(10^{-(n+1)}r,r)$ be the number of balls in a ball covering $B_p(r) \subset \bigcup_{i=1}^N B_{p_i} (10^{-(n+1)}r)$. Then
\[
\mbox{cont} (p,r) \leq (n+1) 2^{N(10^{-(n+1)}r,r)} \max_{p' \in B_p(r), j=1,\cdots, n+1} \mbox{cont} (p', 10^{-j}r).
\]
\end{lem}

The second reduction is the {\em compression}, for its definition we refer to \cite{g,c}. The content of a ball is bounded from above by that of its compression.

Combining the two, we first compress the ball until it is incompressible; then reduce to a ball of one-tenth of the size and compress again. The number of steps to go from $B_p(r)$ to a contractible ball is defined to be the {\em rank}, rank $(p,r)$. (For the precise definition, see \cite{g,c}.) Clearly, by Lemma~\ref{cs} together with the fact that content does not decrease under compression and the defintion of rank, we have \cite[Corollary 5.13]{c}
\begin{lem} \label{cbr}
Let $N(10^{-(n+1)}r,r)$ be as in Lemma~\ref{cs}. Then
\[
\mbox{cont} (p,r) \leq ((n+1) 2^{N(10^{-(n+1)}r,r)})^{rank (p,r)}.
\]
\end{lem}

The number $N(10^{-(n+1)}r,r)$ is known to be uniformly bounded in the presence of a Ricci curvature lower bound.
\begin{prop}[Gromov] (see \cite[Proposition 3.11]{c}) \label{bc}
Let the Ricci curvature of $M^n$ satisfy $\Ric_{M^n} \geq (n-1)H$. Then given $r,\epsilon > 0$ and $p \in M^n$, there exists a covering, $B_p(r) \subset \cup_1^N B_{p_i} (\epsilon), \ (p_i$ in $B_p(r))$ with $N \leq N_1 (n, Hr^2, r/\epsilon)$. Moreover, the multiplicity of this covering is at most $N_2(n,Hr^2)$.
\end{prop}

The key ingredient in the estimate of the rank is the following lemma of 
 Gromov. In the Gromov's proof of the uniform betti number estimate (Theorem~\ref{gtheo}), the sectional curvature hypothesis is only used here.
\begin{lem}  [Gromov's lemma]
 Let $q_1$ be critical with respect to $p$ and let $q_2$ satisfy
\[
|pq_2| \geq \nu |pq_1|,
\]
for some $\nu > 1$. Let $\gamma_1, \gamma_2$ be minimal geodesics from $p$ to $q_1,q_2$ respectively and put $\theta = \angle (\gamma_1'(0), \gamma_2'(0))$. If $K_M \geq H, (H < 0)$ and $|pq_2| \leq D$, then
\[ \theta \geq \cos^{-1} (\frac{\tanh (\sqrt{-H} D/\nu )}{\tanh (\sqrt{-H}D)}).
\]
\end{lem}

Correspondingly we prove the following local version of this lemma for Ricci curvature.
\begin{lem} \label{cpe}
 Let $q_1$ be critical with respect to $p$ and let $q_2$ satisfy
\[
|pq_2| \geq \nu |pq_1|,
\]
for some $\nu > \frac{1+\sin \frac{\pi}{36}}{1-\sin \frac{\pi}{36}} = 1.190954\cdots$. Let $\gamma_1, \gamma_2$ be minimal geodesics from $p$ to $q_1,q_2$ respectively and put $\theta = \angle (\gamma_1'(0), \gamma_2'(0))$. If $\Ric_M \geq (n-1)H$, conj $\geq r_0$ and $|pq_2| \leq r_{ac} (p)/4$, then
\[ \theta \geq \frac{18}{19}\cos^{-1} (\sin \frac{\pi}{36} + (1+ \sin \frac{\pi}{36} ) \frac{1}{\nu}).  \]
\end{lem}
\Pf Put $|pq_1| = x, \ |q_1q_2| = y,\ |pq_2| = z$. Let $\sigma$ be minimal from $q_1$ to $q_2$.  Since $q_1$ is critical to $p$, there exists $\tau$, minimal from $q_1$ to $p$ with 
\[
\theta_1 = \angle (\sigma'(0), \tau'(0)) \leq \frac{\pi}{2}.
\]
Let $\bar{\theta},\ \bar{\theta}_1$ be the corresponding angles in $R^2$. Then by applying
Theorem~\ref{tce} to the geodesic triangle $\{ \tau, \gamma_2, \sigma \}$, we have
\[
\theta_1 \geq \frac{18}{19} \bar{\theta}_1.
\]
Therefore
\[
\bar{\theta}_1 \leq \frac{19}{36} \pi.
\]
Now
\ba     \label{l1}
z^2  & =&  x^2 + y^2 -2xy \cos \bar{\theta}_1  \nonumber \\
& \leq & x^2 + y^2 -2xy \cos (\frac{19\pi}{36}),
\ea
and
\be \label{l2}
y^2  =  x^2 + z^2 -2xz \cos \bar{\theta}.
\ee
Combining (\ref{l1}) and (\ref{l2}) gives
\[
z\cos \bar{\theta} \leq x - y \cos \frac{19\pi}{36},
\]
which gives
\ban
\cos \bar{\theta} & \leq & \frac{x + y \sin \frac {\pi}{36}}{z} \\
& \leq &   \frac{x + (x+z) \sin \frac {\pi}{36}}{z} \ \ \ \ \mbox{since}\  \ y \leq x+z \\
& \leq &  \sin \frac {\pi}{36} + (1+  \sin \frac {\pi}{36} ) \frac{1}{\nu} \ \ \ \ \mbox{since}\ \  \frac{x}{z} \leq \frac{1}{\nu}.
\ean
Therefore
\[
\bar{\theta} \geq \cos^{-1} ( \sin \frac {\pi}{36} + (1+  \sin \frac {\pi}{36} ) \frac{1}{\nu}).
\]
Now applying Theorem~\ref{tce} to the geodesic triangle $\{\gamma_1,\gamma_2,\sigma \}$ gives
\[
\theta \geq \frac{18}{19} \bar{\theta} \geq \frac{18}{19} \cos^{-1} ( \sin \frac {\pi}{36} + (1+  \sin \frac {\pi}{36} ) \frac{1}{\nu}).
\]
\qed
\begin{coro}  \label{nume}
Let $q_1,\cdots,q_N$ be a sequence of critical points of $p$, with 
\[
|pq_{i+1}| \geq \mu |pq_i| \ \ (\mu >1.191).
\]
If $\Ric_M \geq H$, conj $\geq r_0$ and $|pq_N| \leq r_{ac}(p)/4$, then
\[
N \leq {\cal N} (n,\mu,Hr_{ac}^2).
\]
\end{coro}
\Pf This is standard. We follow the one in \cite{c}. Take minimal geodesics, $\gamma_i$ from $p$ to $q_i$. View $\{\gamma_i'(0)\}$ as a subset of $S^{n-1} \subset M_p^n$. Then Lemma~\ref{cpe} gives a lower bound on the distance, $\theta$, between any pair $\gamma_i'(0), \gamma_j'(0)$. The balls of radius $\theta/2$ about $\gamma_i'(0) \in S^{n-1}$ are mutually disjoint. Hence, if we denote by $V_{n-1,1}(r)$, the volume of a ball of radius $r$ on $S^{n-1}$, we can take
\[
{\cal N} = \frac{V_{n-1,1}(\pi)}{V_{n-1,1}(\theta/2)},
\]
where $V_{n-1,1}(\pi) =$ Vol$ (S^{n-1})$ and $\theta$ is the minimum value allowed by Lemma~\ref{cpe}.
\qed

Now the size of rank$(p,r)$ is related to the existence of critical points as follows \cite[Lemma 6.4]{c}.
\begin{lem}  \label{c-r}
Let $M^n$ be Riemannian and let rank $(r,p) = j$. Then there exists $y \in B_p(5r)$ and $x_j, \cdots, x_1 \in B_p(5r)$, such that for all $i \leq j,\ x_i$ is critical with respect to $y$ and
\[
|x_iy| \geq \frac{5}{4} |x_{i-1}y|.
\]
\end{lem}

Combining Lemma~\ref{c-r}, Corollary~\ref{nume} and Lemma~\ref{cbr} gives
\begin{prop}  \label{lc}
For manifolds satisfying (\ref{class}) 
\[
\mbox{cont} (p,r) \leq C(n, r_{ac}^2H) \ \ \mbox{for all} \ \ r\leq r_{ac}/4,
\]
where $r_{ac}$ is the angle comparison radius defined in \S 2.
\end{prop}

Using Proposition~\ref{lc} and Lemma~\ref{cs} inductively completes the proof of Theorem~\ref{bn}.

\sect{Uniform estimate of generators of fundamental groups}
As in \cite{g2,bk} we use the short basis trick to prove Theorem~\ref{fg}. To localize the problem, we introduce
\[
\pi_1(p,r) = \mbox{Im}\ (\pi_1 (B_p(r) \ra \pi_1(M)),
\]
where the map is induced by inclusion. The following lemma is a weak version of  Van Kampen theorem which can be found in \cite{m}.
\begin{lem}
If $\{B_{p_i} (r)\}$ is a ball covering of $M$, then $\pi_1(M)$ is generated by $\{ \pi_1 (p_i,r)\}$.
\end{lem}

 From this lemma and Proposition~\ref{bc}, we see that it suffices to establish the estimate for generators of $\pi_1(p,r_1)$ for all $p \in M$ and $r_1 = \frac{1}{6} \min \{ r_{ac}, r_0 \}$. To show this, first we have
\begin{lem}  \label{gfg}
The fundamental group $\pi_1 (B_p(r_1))$ (and hence $\pi_1(p,r_1)$) is generated by geodesic loops at $p$ of length $< 2r_1$.
\end{lem}
\Pf For any closed curve $c$ at $p$ contained in $B_p(r_1)$, we can choose an $r_1' < r_1$ such that $c$ is contained in $B_p(r_1')$. Now the same argument as in the proof of Proposition 2.1.5 of \cite{bk} shows that $c$ is represented as a product of geodesic loops at $p$ of length $\leq 2r_1' < 2r_1$.
\qed

Now represent each element of $\pi_1(p,r_1)$ by a shortest geodesic loop $\alpha$ at $p$ and call $|\alpha| = l(\alpha)$ the length of the homotopy class. We point out that this may not be possible for $\pi_1(B_p(r_1))$. Note also that $\alpha$ is not necessarily contained in $B_p(r_1)$. We now pick a short basis $\{\alpha_1, \cdots, \alpha_s \}$ for $\pi_1(p,r_1)$ as follows:

1) $\alpha_1$ represents a nontrivial homotopy class of minimal length.

2) If $\alpha_1, \cdots, \alpha_k$ have already been chosen, then $\alpha_{k+1}$ represents a homotopy class of {\em minimal} length in the {\em complement} of the subgroup generated by $\{\alpha_1, \cdots, \alpha_k \}$.

\begin{lem}
We have  
\be \label{le}
|\alpha_1| \leq |\alpha_2| \leq \cdots
\ee
 and 
\be  \label{ge}
|\alpha_i \alpha_j^{-1}| \geq \max \{ |\alpha_i|, |\alpha_j| \}.
\ee
Furthermore,
\[
|\alpha_i| < 2r_1.
\]
\end{lem}
\Pf
(\ref{le}) and (\ref{ge}) are clear by definition (otherwise $\alpha_i$ or $\alpha_j$ were not chosen minimally). To see the last statement we assume on the contrary that
\[
|\alpha_k | \geq 2r_1.
\]
By Lemma~\ref{gfg}, $[\alpha_k]$ is a product of homotopy classes represented by geodesic loops of length $< 2r_1$. At least one of the class, say $[\alpha_k']$ is in the complement of the subgroup $\lp \alpha_1, \cdots, \alpha_{k-1} \rp$, otherwise $[\alpha_k ] \in \lp \alpha_1, \cdots, \alpha_{k-1} \rp$ contradicting the choice of $\alpha_k$. But then 
\[
|\alpha_k'| < 2r_1 \leq |\alpha_k|
\]
also contradicting the minimality of $\alpha_k$.
\qed

\begin{lem}
We have uniform bound on the number of elements in a short basis $\{\alpha_1, \cdots, \alpha_s \}$ of $\pi_1 (p, r_1)$:
\be \label{us}
s \leq (\frac{19}{3})^{n-1}.
\ee
\end{lem}
\Pf By definition $\alpha_i, \alpha_j$ are geodesic loops of $M$ at $p$ of length $|\alpha_i|, |\alpha_j|$ respectively. Moreover $\alpha_i\alpha_j^{-1}$ is homotopic to a geodesic loop $\alpha_{ij}$ of $M$ of length $|\alpha_i\alpha_j^{-1}|$. We lift $\alpha_i,\ \alpha_j,\ \alpha_{ij}$ to the universal covering $\tilde{M}$ and obtain a triangle with edge lengths $|\alpha_i|, |\alpha_j|, |\alpha_i \alpha_j^{-1}|$. Note that the triangle is contained in $B_{\tilde{p}} (6r)$ and the conjugate radius of $\tilde{M}$ is also equal to $r_0$. Thus Theorem~\ref{tce} applies. In particular the angle opposite to $|\alpha_i \alpha_j^{-1}|$ is greater than or equal to $\frac{18}{19} \frac{\pi}{3} = \frac{6\pi}{19}$.

Now the proof of Corollary~\ref{nume} gives (\ref{us}).
\qed

\sect{Fibration theorem}
We first construct the map $f:M \ra N$ as in \cite{f,y}. 

By the definition of the Hausdorff distance, if $d_H(M,N) < \epsilon$, then there exsits a metric $d$ on the disjoint union of $M$ and $N$ satisfying the following: \\
1) The restriction of $d$ to $M$ and $N$ coincide with the original Riemannian distance of $M$ and $N$. \\
2) For every $x \in M, y' \in N$, there exsit $x' \in N, y\in M$ such that $d(x,x') < \epsilon, \ d(y,y') < \epsilon$.  \\
Hence we can take a discrete subset $\{m_i\} \subset M$ and $\{m_i'\} \subset N$ such that \\
1) They are $7\epsilon$-dense in $M$ and $N$ respectively, \\
2) $d(m_i,m_j) > \epsilon, \ d(m_i', m_j') > \epsilon$ for all $i \not= j$, \\
3) $d(m_i, m_i') < \epsilon$.

Let $R = \frac{1}{4} \min \{i_0, r_0, r_{ac} \}$, $S = \# \{m_i\}=\#\{m_i'\}$, $\sigma < R$ and $\chi : [0, \infty) \ra [0, \ 1]$ a smooth function such that
\[ \left\{ \begin{array}{ll}
\chi (t) = 0 & \mbox{if} \ \ t \geq \sigma \\
\chi (t) = 1 & \mbox{if} \ \ t \leq \hf \sigma \\
\chi'(t) < 0 & \mbox{if} \ \ \hf \sigma < t < \sigma.
\end{array} \right. \]
Then we define  $f_N: N \ra R^S$ and $f_M: M \ra R^S$ by
 \ban
 f_N(x) & = & (\chi (d(x, m_i')))_i \\
f_M(x) & = & (\chi\left(\frac{\int_{y \in B_{m_i}(\epsilon)} d(x,y) dy}{ \vol (B_{m_i}(\epsilon))}\right))_i.
\ean
Since $|K_N| \leq 1$ and inj$_N \geq i_0$ $f_N$ is an embedding. Let ${\cal N} (N)$ be the normal bundle of $f_N(N)$ in $R^S$ and ${\cal N}_\delta(N) = \{ v\in {\cal N}(N)| |v| < \delta\}$. Then there exsits $\delta > 0$ such that the normal exponential map restricted to ${\cal N}_\delta(N)$ is a diffeomorphism between ${\cal N}_\delta(N)$ and $B_\delta (f_N(N))$, the $\delta$-neighborhood of $f_N(N)$. Let $\pi : B_\delta (f_N(N)) \ra f_N(N)$ be the projection along fibre of normal bundle. Now the map $f = f_N^{-1} \circ \pi \circ f_M: M\ra N$ is well defined. From the construction we see that $d(x, f(x)) < \tau(\epsilon)$ for each $x \in M$. Here $\tau (\epsilon)$  is a positive number depending only on $\epsilon, \ n, \ i_0, \ r_0$ and satisfying $\lim_{\epsilon \ra 0} \tau (\epsilon) = 0$. From now on, we use a notation $\tau (a, \cdots, b |c)$ to denote a positive number depending only on $n, \ i_0, \ r_0$ and $a, \cdots, b,c$ and satisfying $\lim_{c \ra 0}\tau (a, \cdots, b |c) = 0$ for fixed $a, \cdot, b$.

To prove that $f$ is a fibration, it suffices to see that $f$ is of maximal rank. The basic lemma is the following (This corresponds to Lemma 2.6 in \cite{y} and Lemma 2.1 in \cite{f})
\begin{lem} \label{btheta}
Let $\nu, l, l'$ be positive numbers with $l \leq l' <R$, and let $c_i: [0,t_i] \ra M$ and $c_i': [0, t_i'] \ra N (i =1,2)$ be minimal geodesics with $c_1(0) = c_2(0),\ c_1'(0) = c_2'(0)$ such that
\[
d(c_i(0), c_i'(0)) < \nu, \ d(c_i(t_i), c_i'(t_i)) < \nu, \ l/10 < t_1 < l, \ l'/10 < t_2 < l'.
\]
Denote $\theta = \angle (\frac{Dc_1}{dt} (0), \frac{Dc_2}{dt} (0)), \theta' = \angle (\frac{Dc_1'}{dt} (0), \frac{Dc_2'}{dt} (0))$.
Then we have
\be \label{theta}
|\theta - \mu \theta'| \leq (1-\mu) \pi + \tau (l) + \tau (l,l'|\nu) + \tau (l,l'|\epsilon).
\ee
\end{lem}
\Pf As in \cite[Lemma 2.5]{y} we may assume $l = l'$. Let $c$ be the geodesic such that $c(0) = c_1(0), \ \frac{Dc}{dt}(0) = -\frac{Dc_1}{dt}(0)$. Take a point $x' \in N$ such that $d(x', c(t_1)) < \epsilon$. Let $c': [0, t_0'] \ra N$ be a minimal geodesic joining $c_i'(0)$ to $x'$. Then by Theorem~\ref{tce}
\[
\theta \leq \mu \bar{\theta},
\]
where $\bar{\theta}$ denote the comparison angle in the Euclidean plane. By the Toponogov comparison theorem
\[
\theta' \leq \bar{\theta'},
\]
where $\bar{\theta'}$ denote the comparison angle in the sphere.
Therefore
\be  \label{theta1}
\theta \leq  \mu \theta' - \tau(l) -\tau(l|\nu).
\ee
Denote $\theta_1 =  \angle (\frac{Dc}{dt} (0), \frac{Dc_2}{dt} (0)),\ \theta_1'
=\angle (\frac{Dc'}{dt} (0), \frac{Dc_2'}{dt} (0))$. Then similarly we have
\begin{equation}
\label{theta2}
\theta_1 \leq  \mu \theta_1' - \tau(l) -\tau(l|\nu) - \tau (l|\epsilon).
\end{equation}
Since $|d(c'(t_0'),c_1'(t_1')) - (t_0'+t_1')| < 2\epsilon + 4 \nu$, we have
\be  \label{theta3}
|\angle (\frac{Dc_1'}{dt}(0), \frac{Dc'}{dt}(0)) -\pi | \leq \tau(l|\epsilon) + \tau(l |\nu).
\ee
Note that $\theta + \theta_1 = \pi, \ \angle (\frac{Dc_1'}{dt}(0), \frac{Dc'}{dt}(0)) \leq \theta' + \theta_1'$. Combining inequalities (\ref{theta1}), (\ref{theta2}), (\ref{theta3}) yields (\ref{theta}).
\qed

By Lemma~\ref{btheta} and the construction of $f_N$ and $f_M$, we have the following lemma (for details see \cite[Lemma 2.8]{y}).
\begin{lem}
For each $x \in M$ and unit vector $\xi' \in T_{f(x)}N$, let $c': [0, t'] \ra N$ be the geodesic such that $\frac{Dc'}{dt}(0) = \xi', \ l < t' < R, \ l \geq \sigma$. Take a minimal geodesic $c: [0,t] \ra M$ such that $c(0) = x, \ d(c(t), c'(t')) < \epsilon$ and put $\xi = \frac{Dc}{dt}(0)$. Then we have
\[
|df(\xi) - \mu \xi'| < 3(1-\mu)\pi  + \tau (\sigma) + \tau (\sigma, l|\epsilon).
\]
\end{lem}

Now it is clear that when $\mu$ is sufficiently close to $1$ and $\sigma, \epsilon$ sufficiently small, $f$ is a submersion, hence a fibration.

\sect{Proof of Proposition~\protect\ref{sph}}
For $p_0,\ p_1 \in M$, the excess function $e_{p_0,p_1}: M^n \ra R$ is defined by \[
e_{p_0,p_1}(p) = d(p_0,p) + d(p_1,p) - d(p_0,p_1).
\]
It measures the ``excess" in the triangle inequality.

Note that the excess function is ``monotonic" in the following sense. If $p_0',\ p_1'$ are any two points lying on a minimal geodesic connecting $p_0,\ p$ and $p_1,\ p$ respectively, then
\be \label{e}
e_{p_0,p_1} (p) \geq e_{p'_0,p'_1} (p).
\ee
This is an easy consequence of the triangle inequality.

The proof of Proposition~\ref{sph} follows the line of \cite{gp}, 
the case with sectional curvature lower bound. First by a volume comparison argument one has
(cf. \cite[Lemma 1]{gp})

\begin{lem}
Let $M$ be a complete Riemannain manifold with $\Ric_M \geq n-1, \ \diam \geq \pi - \epsilon$, and let $p,q \in M$ with $|pq| = \diam_M$, then
\be \label{se}
e_{p,q} (x) \leq \kappa (\epsilon) \ \mbox{for all} \ x \in M,
\ee
where $\kappa (\epsilon)$ is a positive function of $\epsilon$ and $\kappa (\epsilon) \ra 0$ as $\epsilon \ra 0$.
\end{lem}

\noindent {\em Proof of Proposition~\protect\ref{sph}} It suffices to prove that the angle between any minimal geodesic $\gamma_1, \gamma_2$ from $x$ to $p$ and from $x$ to $q$ respectively is $> \frac{\pi}{2}$. Let $p',q' \in \gamma_1, \gamma_2$ such that $|xp'|= |xq'| =\min \{ r_{ac}/4, \delta\}$.
By (\ref{e}) and (\ref{se})
\begin{equation} 
e_{p',q'}(x) \leq e_{p,q}(x) \leq \kappa (\epsilon).   \label{ne}
\end{equation}
Applying Theorem~\ref{tce} to the triangle $p'xg'$, we have
\begin{eqnarray*}
\angle p'xq' & \geq & \frac{18}{19} \cos^{-1} \frac{|xp'|^2 +|xq'|^2-|p'q'|^2}{ 2|xp'||xq'|} \\
& = & \frac{18}{19}\cos^{-1} (-1- \frac{e^2_{p',q'}(x) -4 e_{p',q'}(x) |xp'|}{2|xp'|^2}) \\
& \geq &  \frac{18}{19}\cos^{-1} (-1 + \frac{2\kappa (\epsilon)}{|xp'|}).
\end{eqnarray*}
Here the last inequality follows from (\ref{ne}).
Now it is clear that there is an $\epsilon (\delta, r_0)$ such that the angle $\angle pxq = \angle p'xq' > \frac{\pi}{2}$.

\qed

\end{document}